\documentclass[a4paper,12pt]{elsarticle}
\usepackage{amsmath}
\usepackage[english]{babel}
\usepackage{graphicx}
\usepackage[usenames,dvipsnames]{color} 
\usepackage[numbers]{natbib}
\usepackage[footnotesize,bf,sf,center]{caption}
\usepackage{float}
\usepackage[dvips]{epsfig}
\usepackage[bottom=2.5cm,top=2.5cm,left=3cm,right=2cm]{geometry}
\usepackage{calc}
\usepackage{soul}
\usepackage{amsmath,amsfonts,amssymb}
\usepackage{subfigure}
\usepackage{setspace}
\usepackage{amssymb}
\usepackage{mathrsfs}
\usepackage{algorithm}
\usepackage{algpseudocode}
\usepackage{dsfont}
\usepackage{epstopdf}
\usepackage{subfig}%
\usepackage[utf8]{inputenc}
\usepackage{amsfonts,amsmath}
\newcommand{\ve}[1]{\mbox{\boldmath $#1$}}
\newtheorem{theorem}{Theorem}[section]
\newdefinition{rmk}{Remark}

\newtheorem{lemma}[theorem]{Lemma}

\newcommand{\proof} [1]{ \noindent {\bf Proof.} #1 \hfill\rule{0.5em}{1.2ex} \par\medskip}
\usepackage{graphicx}
\newcommand{\vvvert}[1]{{\left\vert\kern-0.25ex\left\vert\kern-0.25ex\left\vert #1 
    \right\vert\kern-0.25ex\right\vert\kern-0.25ex\right\vert}}
\begin{document}
\newcommand{\red}[1]{\textcolor{red}{#1}}
\newcommand{\blue}[1]{\textcolor{blue}{#1}}
\newcommand{\Orient}[1]{\colorbox{yellow}{\bf \ ... }\textcolor{ForestGreen}{#1}\colorbox{yellow}{\bf ... \ }}
\newcommand{\vol}{\mathop{\ooalign{\hfil$V$\hfil\cr\kern0.08em--\hfil\cr}}\nolimits}
\begin{frontmatter}

  \renewcommand\arraystretch{1.0}

    \title{\textbf{Fully consistent lowest-order finite element methods for generalised Stokes flows with variable viscosity}}
    


    \author{
   {\bf Felipe Galarce} $^{1}$, \ 
   {\bf Douglas R.~Q.~Pacheco} $^{2,3,4}$
   \footnote[0]{{\sf Email address:} {\tt pacheco@ssd.rwth-aachen.de}, 
     {\sf corresponding author}}\\
      {\small ${}^{1}$ School of Civil Engineering, Pontificia Universidad Católica de Valparaíso, Chile}\\
   {\small ${}^{2}$ Chair for Computational Analysis of Technical Systems, RWTH Aachen University, Germany}\\
   {\small ${}^{3}$ Chair of Methods for Model-based Development in Computational Engineering, RWTH}\\
   {\small ${}^{4}$ Center for Simulation and Data Science (JARA-CSD), RWTH Aachen University, Germany}
    }


\begin{abstract}
Variable viscosity arises in many flow scenarios, often imposing numerical challenges. Yet, discretisation methods designed specifically for non-constant viscosity are few, and their analysis is even scarcer. In finite element methods for incompressible flows, the most popular approach to allow equal-order velocity-pressure interpolation are residual-based stabilisations. For low-order elements, however, the viscous part of that residual cannot be approximated, often compromising accuracy. Assuming slightly more regularity on the viscosity field, we can construct stabilisation methods that fully approximate the residual, regardless of the polynomial order of the finite element spaces. This work analyses two variants of this fully consistent approach, with the generalised Stokes system as a model problem. We prove unique solvability and derive expressions for the stabilisation parameter, generalising some classical results for constant viscosity. Numerical results illustrate how our method completely eliminates the spurious pressure boundary layers typically induced by low-order PSPG-like stabilisations. 
\end{abstract}
\end{frontmatter}

\section{Introduction}
Variable viscosity fields appear in various flow scenarios, such as when non-Newtonian behaviour, turbulence models or temperature gradients are considered. Nonetheless, in most of the Computational Fluid Dynamics literature, the assumption of constant viscosity is often taken for granted. Yet, it is not always straightforward to extend classical numerical techniques to fluid problems with non-constant viscosity. Even a simple, quasi-linear rheological law can lead to numerical instabilities, spurious boundary conditions or undesired matrix structures \cite{Pacheco2021CMAME}. In that context, the literature on numerical methods for variable-viscosity flow problems has recently been expanding. Some examples are stabilisation methods \cite{Schussnig2021,Barrenechea2023}, fractional-step schemes \cite{Pacheco2021CMAME,Deteix2018,Plasman2020,Barrenechea2024} and mixed formulations \cite{Anaya2021,Anaya2023}.

For incompressible flows with variable viscosity, we have recently presented \cite{Schussnig2021} an equal-order stabilised method that solves an old issue of residual-based methods: the weakened consistency and resulting numerical boundary layers observed when low-order elements are used \cite{Jansen1999}. When the polynomial order of the finite element spaces is smaller than the highest derivative in the equations, part of the stabilising residual vanishes, causing spurious pressure boundary layers \cite{Burman2006}. Our modified approach, on the other hand, assumes Lipschitz continuity on the viscosity field to rewrite the viscous term using only first-order derivatives. The resulting formulation has a similar structure and cost as standard residual-based methods, but is fully consistent even for piecewise linear elements. Despite previous numerical experiments showcasing the accuracy and stability of the method \cite{Schussnig2021,Pacheco2021}, an a-priori analysis has not been presented yet, which is thus the main goal of this work. We prove the unique solvability of the modified stabilisation, derive closed-form expressions for the stabilisation parameter and sketch convergence. Numerical examples showcase that our method indeed solves the issue of spurious boundary layers and loss of pressure accuracy on the boundaries.

\section{Preliminaries}
Let us consider a domain $\Omega\subset\mathbb{R}^{d}$, $d=2$ or 3, with Lipschitz boundary $\Gamma=\partial\Omega$. As a model problem, we take the generalised Stokes system:
\begin{align}
\sigma\ve{u} - \nabla\cdot(2\nu\nabla^{\mathrm{s}}\ve{u}) + \nabla p &= \ve{f} && \text{in} \ \ \Omega\, ,\label{momentum}\\
\nabla\cdot\ve{u} &= 0  && \text{in} \ \ \Omega\, , \label{incompressibility}\\
\ve{u} &= \mathbf{0} && \text{on} \ \ \Gamma\, ,\label{DirichletBC0}
\end{align}
where the unknowns are the velocity vector $\ve{u}$ and the pressure $p$, while the remaining data are a force vector $\ve{f}$, a viscosity field $\nu$ and a positive constant $\sigma$. The variable viscosity $\nu$ is assumed to satisfy
\begin{equation}
    0 <  \nu_{\mathrm{min}} \leq \nu \leq \nu_{\mathrm{max}} < \infty\, ,
    \label{viscosityBounds}
\end{equation}
and the symmetric gradient operator is defined in the usual way: $2\nabla^{\mathrm{s}} := \nabla + \nabla^{\top}$. 

We consider usual notation for Hilbert and Lebesgue spaces, and denote as $L^2_0(\Omega)$ the space of $L^2(\Omega)$ functions with zero mean on $\Omega$.
The notation $(\cdot\,  ,\cdot)_{\omega}$ will be used for the $L^2(\omega)$ inner product, omitting the domain dependency if $\omega = \Omega$. By $\| \cdot \|$ and $\| \cdot \|_{\infty}$ we denote, respectively, the $L^2(\Omega)$ and $L^{\infty}(\bar{\Omega})$ norms.

Let us assume shape-regular, globally quasi-uniform triangulations $\mathcal{T}$ of $\Omega$, with  
\begin{align}
1 \leq \frac{h}{h_e} \leq Q, \quad Q \geq 1\, ,\label{quasiUniform}
\end{align}
where $h_e$ denotes the size of element $\Omega_e$, and $h$ is the maximum element size. Based on that, we consider standard $H^1$-conforming finite element spaces. The discrete velocity and pressure are sought, respectively, in the finite-element spaces $X_h\subset [H^1_0(\Omega)]^d$ and $Y_h\subset L^2_0(\Omega)$.

\section{Residual-based finite element stabilisations}
The Stokes system being a saddle-point problem, the stability of its discretisation is subject to an inf-sup condition \cite{John2016}. Although it is well known that equal-order velocity-pressure pairs violate that condition, such spaces are attractive due to their simplicity. Stabilisation methods are thus often used to circumvent the discrete inf-sup condition.

\subsection{The classical PSPG method}
The most popular stabilised formulation for equal-order elements is the pressure-stabilised Petrov--Galerkin (PSPG) method: find $\ve{u}_h\in X_h$ and $p_h\in Y_h$ such that 
\begin{align}
\sigma(\ve{u}_h,\ve{v}_h)+\left(2\nu\nabla^{\mathrm{s}}\ve{u}_h,\nabla^{\mathrm{s}}\ve{v}_h\right)  -(p_h,\nabla\cdot\ve{v}_h) &=  (\ve{v}_h,\ve{f}) \, , \label{weakMomentumSD}\\
(q_h,\nabla\cdot\ve{u}_h) + \sum_{\Omega_e\subset\mathcal{T}}(\delta\nabla q_h,\nabla p_h + \sigma\ve{u}_h - \nabla\cdot(2\nu\nabla^{\mathrm{s}}\ve{u}_h)-\ve{f})_{\Omega_e} &= 0
\label{PSPG}
\end{align}
holds for all $\ve{v}_h\in X_h$ and $q_h\in Y_h$, where $\delta$ is a mesh-dependent stabilisation parameter. Although the stabilisation term is proportional to the residual of the momentum equation, a part of the viscous term is lost when considering piecewise linear finite element spaces. Within each element $\Omega_e$, we have
\begin{align*}
    \nabla\cdot(2\nu\nabla^{\mathrm{s}}\ve{u}_h) &\equiv 2\nabla^{\mathrm{s}}\ve{u}_h\nabla\nu + \nu[\Delta\ve{u}_h + \nabla(\nabla\cdot\ve{u}_h)]\\
    &= 2\nabla^{\mathrm{s}}\ve{u}_h\nabla\nu\, . 
\end{align*}
That is, the Laplacian vanishes for the discrete velocity $\ve{u}_h$, which is generally not true for the continuous solution $\ve{u}$. This well-known issue leads to artificial pressure boundary layers and loss of accuracy, especially for coarse meshes \cite{Burman2006}. Yet, we will show how it is possible to overcome that by assuming a higher regularity on the viscosity field.

\subsection{The boundary vorticity stabilisation}
Assuming a sufficiently regular viscosity field, it is possible to use the incompressibility constraint \eqref{incompressibility} to rewrite the stress divergence as
\begin{align*}
\nabla\cdot\left(2\nu\nabla^{\mathrm{s}}\ve{u}\right) &\equiv  2\nabla^{\mathrm{s}}\ve{u}\nabla\nu + \nu\left[2\nabla(\nabla\cdot\ve{u}) - \nabla\times(\nabla\times\ve{u})\right] \\
&= 2\nabla^{\mathrm{s}}\ve{u}\nabla\nu - \nu\nabla\times(\nabla\times\ve{u}) \, . 
\end{align*}
We can use this to construct a stabilised method that preserves full consistency even for lowest-order velocity approximations. For that, a key point is to use the residual not only in element interiors, but in the entire domain $\Omega$, which will require $\ve{f}\in [L^2(\Omega)]^d$ and $\nu\in W^{1,\infty}(\Omega)$. Denoting by $\ve{n}$ the outward unit normal vector on $\Gamma$, we can write
\begin{equation}
\begin{split}
(\nabla q_h,\nabla\cdot(2\nu\nabla^{\mathrm{s}}\ve{u})) &= 
(\nabla q_h,2\nabla^{\mathrm{s}}\ve{u}\nabla\nu) - (\nabla q_h,\nu\nabla\times(\nabla\times\ve{u})) \\
&= (\nabla q_h,2\nabla^{\mathrm{s}}\ve{u}\nabla\nu)+(\ve{n}\times\nabla q_h,\nu\nabla\times\ve{u})_{\Gamma} -(\nabla\times(\nu\nabla q_h),\nabla\times\ve{u}) \\
&= (\nabla q_h,2\nabla^{\mathrm{s}}\ve{u}\nabla\nu)+(\ve{n}\times\nabla q_h,\nu\nabla\times\ve{u})_{\Gamma} +(\nabla q_h\times\nabla\nu,\nabla\times\ve{u})\\
&\equiv (\nabla q_h,2\nabla^{\mathrm{s}}\ve{u}\nabla\nu)+(\ve{n}\times\nabla q_h,\nu\nabla\times\ve{u})_{\Gamma} +(\nabla q_h,\nabla\nu\times(\nabla\times\ve{u}))\\
&\equiv (\nabla q_h,2\nabla^{\mathrm{s}}\ve{u}\nabla\nu)+(\ve{n}\times\nabla q_h,\nu\nabla\times\ve{u})_{\Gamma} +(\nabla q_h,(\nabla^{\top}\ve{u}-\nabla\ve{u})\nabla\nu)\\
&= (\nabla q_h,2\nabla^{\top}\ve{u}\nabla\nu)+(\ve{n}\times\nabla q_h,\nu\nabla\times\ve{u})_{\Gamma}\, ,
\label{transformations}
\end{split}
\end{equation}
where we have used the following identities:
\begin{align*}
    \nabla\times \nabla q_h &\equiv \ve{0} \quad \text{for any} \ \, q_h\in Y_h\subset H^1(\Omega)\, ,\\
    (\ve{a}\times\ve{b})\cdot\ve{c} &\equiv \ve{a}\cdot(\ve{b}\times\ve{c})\quad \text{for any} \ \, \ve{a},\ve{b},\ve{c}\in \mathbb{R}^d\, ,\\
    \ve{a}\times(\nabla\times\ve{b}) &\equiv (\nabla^{\top}\ve{b}-\nabla\ve{b})\ve{a}\quad \text{for any (regular enough)} \ \, \ve{a},\ve{b}\in \mathbb{R}^d\, .
\end{align*}
Thus, by using integration by parts and appropriate identities, we have eliminated the second-order derivatives in a \textsl{consistent} way, so that the final expression in \eqref{transformations} can be fully approximated by the finite element solution $\ve{u}_h$. We can now replace the stabilised equation \eqref{PSPG} by
\begin{equation}
(q_h,\nabla\cdot\ve{u}_h) + \delta(\nabla q_h,\nabla p_h + \sigma\ve{u}_h - 2\nabla^{\top}\ve{u}_h\nabla\nu-\ve{f}) + \delta(\nabla q_h\times\ve{n},\nu\nabla\times\ve{u}_h)_{\Gamma} = 0\, .
\label{BVS}
\end{equation}
We denote this as the boundary vorticity stabilisation (BVS) method, due to the integral on $\Gamma$ depending on $\nabla\times\ve{u}_h$. Notice that all the integrals are well defined when using standard continuous finite element spaces for $(\ve{v}_h,\ve{u}_h,q_h,p_h)$.

\section{Stability analysis}
\subsection{Stress-divergence formulation}
The first formulation we will analyse is the stress-divergence (SD) variant, that is, combining the momentum equation \eqref{weakMomentumSD} with the BVS \eqref{BVS}. 
\begin{lemma}[Unique solvability of the BVS formulation in SD form]
For \\ $\ve{f}\in [L^2(\Omega)]^d$ and $\nu\in W^{1,\infty}(\Omega)$ bounded as in \eqref{viscosityBounds}, problem \eqref{weakMomentumSD},\eqref{BVS} has a unique solution if the stabilisation parameter $\delta$ satisfies 
\begin{align}
    \delta \leq \mathrm{min}\left\lbrace\frac{1}{3\sigma}, \, \frac{\nu_{\mathrm{min}} h^2/12}{h^2\|\nabla\nu\|^2_{\infty} + C\nu_{\mathrm{max}}^2} \right\rbrace , 
    \label{deltaSD}
\end{align}
with $C$ being a positive constant depending only on the geometry and the discretisation.
\end{lemma}
\proof{To prove the unique solvability, we first show that the corresponding bilinear form 
\begin{align*}
  &A_{\mathrm{SD}}((\ve{v}_h,q_h),(\ve{u}_h,p_h)) := \sigma(\ve{u}_h,\ve{v}_h)+\left(2\nu\nabla^{\mathrm{s}}\ve{u}_h,\nabla^{\mathrm{s}}\ve{v}_h\right)  -\left(p_h,\nabla\cdot\ve{v}_h\right)  \\  &+(q_h,\nabla\cdot\ve{u}_h) + \delta(\nabla q_h,\nabla p_h + \sigma\ve{u}_h - 2\nabla^{\top}\ve{u}_h\nabla\nu) + \delta(\nabla q_h\times\ve{n},\nu\nabla\times\ve{u}_h)_{\Gamma}
\end{align*}
is coercive with respect to the norm 
\begin{equation}
    \vvvert{(\ve{v}_h,q_h)} := \left(\sigma\| \ve{v}_h\|^2 + \nu_{\mathrm{min}}\|\nabla\ve{v}_h\|^2 + \delta\|\nabla q_h\|^2\right)^{1/2}\, .
    \label{norm}
\end{equation}
For any $\ve{v}_h\in X_h$ and $q_h\in Y_h$, we have
\begin{align*}
A_{\mathrm{SD}}((\ve{v}_h,q_h),(\ve{v}_h,q_h)) = \  &\sigma\|\ve{v}_h\|^2 + 2\|\nu^{\frac{1}{2}}\nabla^{\mathrm{s}}\ve{v}_h\|^2 + \delta\|\nabla q_h\|^2  \\
&+\delta\big[(\nabla q_h,\sigma\ve{v}_h-2\nabla^{\top}\ve{v}_h\nabla\nu) + (\nabla q_h\times\ve{n},\nu\nabla\times\ve{v}_h)_{\Gamma} \big],
\end{align*}
where some terms must be estimated. Using the boundedness of $\nu$ and Korn's inequality (\citet{John2016}, Lemma 3.37), we can write
\begin{align}
    2\|\nu^{\frac{1}{2}}\nabla^{\mathrm{s}}\ve{v}_h\|^2 \geq 
    2\nu_{\mathrm{min}}\|\nabla^{\mathrm{s}}\ve{v}_h\|^2 \geq
    \nu_{\mathrm{min}}\|\nabla\ve{v}_h\|^2\, .
    \label{Korn}
\end{align}
The Cauchy-Schwarz and Young inequalities lead to
\begin{equation}
\begin{split}
\left|(\nabla q_h,\sigma\ve{v}_h)\right| &\leq \|\nabla q_h\|\sigma\|\ve{v}_h\|\\
&\leq  \frac{1}{6}\|\nabla q_h\|^2 + \frac{3\sigma}{2}\sigma\|\ve{v}_h\|^2
\end{split}
\end{equation}
and
\begin{equation}
\begin{split}
   \left|(\nabla q_h,2\nabla^{\top}\ve{v}_h\nabla\nu)\right| &\leq 2\|\nabla q_h\| \|\nabla\ve{v}_h\|\|\nabla\nu\|_{\infty}\\
   &\leq   \frac{1}{6}\|\nabla q_h\|^2 + 6\|\nabla\nu\|_{\infty}^2\|\nabla\ve{v}_h\|^2 \, .
\end{split}
\end{equation}

A major distinction between the PSPG and BVS methods is that the latter contains a boundary integral. Estimating that requires a trace inequality (\citet{DiPietro2012}, Lemma 1.46):
\begin{align*}
\sqrt{h_e}\| v_h \|_{L^2(F_e)} \leq C_{\text{tr}}\| v_h \|_{L^2(\Omega_e)} \, \ \text{for any face $F_e$ of $\Omega_e$,}
\end{align*}
with $C_{\text{tr}}$ depending only on shape-regularity constants, the spatial dimension $d$ and the polynomial degree of $X_h$. Hence:
\begin{align*}
(\nabla q_h\times\ve{n}, \nu\nabla\times\ve{v}_h)_{\Gamma} &= \sum_{\Omega_e\subset\mathcal{T}}(\nabla q_h \times\ve{n}, \nu\nabla\times\ve{v}_h)_{\Gamma\cap\partial\Omega_e}\\
 &\leq  \sum_{\Omega_e\subset\mathcal{T}}\| \nabla q_h \times\ve{n}\|_{L^2(\Gamma\cap\partial\Omega_e)}\|\nu\nabla\times\ve{v}_h \|_{L^2(\Gamma\cap\partial\Omega_e)} \\
&\leq \nu_{\mathrm{max}}\sum_{\Omega_e\subset\mathcal{T}}\| \nabla q_h \|_{L^2(\Gamma\cap\partial\Omega_e)} \| \nabla \ve{v}_h \|_{L^2(\Gamma\cap\partial\Omega_e)} \\
 &\leq C_{\text{tr}}^2\nu_{\mathrm{max}}\sum_{\Omega_e\subset\mathcal{T}}{h_e}^{-1}\| \nabla q_h \|_{L^2(\Omega_e)} \| \nabla \ve{v}_h \|_{L^2(\Omega_e)} \, .
\end{align*}
Now, using \eqref{quasiUniform} with the triangle and Young inequalities yields
\begin{equation}
\begin{split}
|(\nabla q_h\times\ve{n}, \nu\nabla\times\ve{v}_h)_{\Gamma}| &\leq   QC_{\text{tr}}^2 h^{-1}\nu_{\mathrm{max}}\| \nabla q_h \| \, \| \nabla \ve{v}_h \| \\
&\leq \frac{1}{6}\|\nabla q_h \|^2 + 6C\frac{\nu_{\mathrm{max}}^2} {h^2}\|\nabla\ve{v}_h\|^2 \, ,\label{term4}
\end{split}
\end{equation}
with $C= C_{\text{tr}}^2Q/2$. Adding up estimates \eqref{Korn}--\eqref{term4} gives us
\begin{align*}
&A_{\mathrm{SD}}((\ve{v}_h,q_h),(\ve{v}_h,q_h)) \geq \frac{1}{2}\vvvert{(\ve{v}_h,q_h)}^2 \\
&+ (1-3\delta\sigma)\frac{\sigma}{2}\|\ve{v}_h\|^2 +  \bigg[1- \frac{12\delta}{h^2\nu_{\mathrm{min}}}(h^2\|\nabla\nu\|^2_{\infty} + C\nu_{\mathrm{max}}^2)\bigg]\frac{\nu_{\mathrm{min}}}{2}\|\nabla\ve{v}_h\|^2\, .
\end{align*}
Thus, by setting the stabilisation parameter according to \eqref{deltaSD}, we have that 
\begin{align}
    A_{\mathrm{SD}}((\ve{v}_h,q_h),(\ve{v}_h,q_h)) \geq \frac{1}{2}\vvvert{(\ve{v}_h,q_h)}^2 
\end{align}
for all $\ve{v}_h\in X_h$ and $q_h\in Y_h$. Finally, the boundedness of both the bilinear form and the right-hand side follows from standard arguments, which guarantees the unique solvability via the thoerem of Lax-Milgram \cite{John2016}.}

\begin{rmk}
For constant $\nu$, the condition on the stabilisation parameter becomes
\begin{align*}
    \delta \leq \mathrm{min}\left\lbrace\frac{1}{3\sigma}, \, \frac{\hat{C}h^2}{3\nu} \right\rbrace, 
\end{align*}
for some positive constant $\hat{C}$. This is, up to the precise determination of $\hat{C}$, the same result one can usually prove for PSPG (see \citet{John2016}, Theorem 5.34).
\end{rmk}

\subsection{On the reaction term in the residual}
As another advantage of using the residual on the entire domain $\Omega$, we can eliminate from the stabilisation parameter the dependence on the reaction coefficient $\sigma$. Integration by parts yields
\begin{align}
    \delta(\nabla q_h,\sigma\ve{u}) = -\sigma\delta(q_h,\nabla\cdot\ve{u}) = 0\, ,
\end{align}
which means we can \textsl{consistently} drop the reaction term from the residual in \eqref{BVS}. As a result, we get an improved condition on the stabilisation parameter:
\begin{align}
    \delta \leq  \frac{\nu_{\mathrm{min}} h^2/12}{h^2\|\nabla\nu\|^2_{\infty} + C\nu_{\mathrm{max}}^2} \, . 
\end{align}

In the more general case of non-homogeneous boundary conditions $\ve{u}|_{\Gamma} = \ve{g}$, one gets
\begin{align*}
    \delta(\nabla q_h,\sigma\ve{u}) = \delta\sigma(q_h,\ve{g}\cdot\ve{n})_{\Gamma}\, ,
\end{align*}
that is, we can still eliminate the term $\delta(\nabla q_h,\sigma\ve{u})$ from the bilinear form, getting a forcing term on the right-hand side instead. 

\subsection{Generalised Laplacian formulation}
Since the viscosity gradient is assumed herein as bounded, a consistent simplification of the viscous term is possible. Condition \eqref{incompressibility} leads to
\begin{align*}
\nabla\cdot\left(2\nu\nabla^{\mathrm{s}}\ve{u}\right) &\equiv \nabla\cdot(\nu\nabla\ve{u}) + \nu\nabla(\nabla\cdot\ve{u}) + \nabla^{\top}\ve{u}\nabla\nu \\
&= \nabla\cdot(\nu\nabla\ve{u}) + \nabla^{\top}\ve{u}\nabla\nu\\
&:= \tilde{\mathbf{D}} \, ,
\end{align*} 
with $\tilde{\mathbf{D}}$ denoting a ``reduced'' stress divergence. Testing with $\ve{v} \in [H^1_0(\Omega)]^d$ yields
\begin{align} 
-(\tilde{\mathbf{D}},\ve{v}) &= (\nu\nabla\ve{u},\nabla\ve{v}) -  (\ve{v},\nabla^{\top}\ve{u}\nabla\nu) - \langle\ve{v},(\nu\nabla\ve{u})\ve{n}\rangle_{\Gamma}\nonumber \\ 
&= (\nu\nabla\ve{u},\nabla\ve{v}) -  (\ve{v},\nabla^{\top}\ve{u}\nabla\nu) 
\, .
\label{genLaplacian}
\end{align}
Although we only address the purely Dirichlet problem herein, the derivation above shows that this formulation induces normal \textsl{pseudo-tractions} as natural boundary conditions, which is often more appropriate to handle outflow boundaries than the real tractions of the classical SD formulation \cite{Heywood1996}. Moreover, the generalised Laplacian (GL) form \eqref{genLaplacian} can be more stable than the SD variant for time-dependent problems \cite{Barrenechea2024}. Motivated by that, we will next analyse the modified problem of finding $\ve{u}_h\in X_h$ and $p_h\in Y_h$ so that
\begin{equation}
\begin{split}
\sigma(\ve{u}_h,\ve{v}_h)+(\nu\nabla\ve{u}_h,\nabla\ve{v}_h) - (\nabla^{\top}\ve{u}_h\nabla\nu,\ve{v}_h)  -(p_h,\nabla\cdot\ve{v}_h) &=  (\ve{f},\ve{v}_h) \, , \\
(q_h,\nabla\cdot\ve{u}_h) + \delta(\nabla q_h,\nabla p_h  - 2\nabla^{\top}\ve{u}_h\nabla\nu) + \delta(\nabla q_h\times\ve{n},\nu\nabla\times\ve{u}_h)_{\Gamma} &= (\ve{f},\delta\nabla q_h)
\label{BVSGL}    
\end{split}
\end{equation}
holds for all $\ve{v}_h\in X_h$ and $q_h\in Y_h$.

\begin{lemma}[Unique solvability of the BVS formulation in GL form]
Assuming $\ve{f}\in L^2(\Omega)$, $\nu\in W^{1,\infty}(\Omega)$ bounded as in \eqref{viscosityBounds}, and 
\begin{align}
    \sigma \geq 3\frac{\|\nabla\nu\|_{\infty}^2}{\nu_{\mathrm{min}}}\, ,
    \label{restriction}
\end{align}
problem \eqref{BVSGL} has a unique solution if
\begin{align}
    \delta \leq  \frac{\nu_{\mathrm{min}} h^2/12}{h^2\|\nabla\nu\|^2_{\infty} + C\nu_{\mathrm{max}}^2}\, ,
    \label{deltaGL}
\end{align}
with $C$ being a positive constant depending only on the geometry and the discretisation.
\end{lemma}
\proof{We will proceed similarly as before, to show that the bilinear form 
\begin{align*}
  &A_{\mathrm{GL}}((\ve{v}_h,q_h),(\ve{u}_h,p_h)) := \sigma(\ve{u}_h,\ve{v}_h)+(\nu\nabla\ve{u}_h,\nabla\ve{v}_h) - (\nabla^{\top}\ve{u}_h\nabla\nu,\ve{v}_h)-\left(p_h,\nabla\cdot\ve{v}_h\right)  \\  &+(q_h,\nabla\cdot\ve{u}_h) + \delta(\nabla q_h,\nabla p_h - 2\nabla^{\top}\ve{u}_h\nabla\nu) + \delta(\nabla q_h\times\ve{n},\nu\nabla\times\ve{u}_h)_{\Gamma}
\end{align*}
is coercive with respect to the norm \eqref{norm}. We have 
\begin{align*}
A_{\mathrm{GL}}((\ve{v}_h,q_h),(\ve{v}_h,q_h)) = \ & \sigma\|\ve{v}_h\|^2 + \|\nu^{\frac{1}{2}}\nabla\ve{v}_h\|^2 + \delta\|\nabla q_h\|^2 \\  & + \delta\big[(\nabla q_h\times\ve{n},\nu\nabla\times\ve{v}_h)_{\Gamma}-(\nabla q_h,2\nabla^{\top}\ve{v}_h\nabla\nu)\big] - (\nabla^{\top}\ve{v}_h\nabla\nu,\ve{v}_h)\, .
\end{align*}
The terms in the bottom row can be estimated as
\begin{align*}
   &\left|(\nabla^{\top}\ve{v}_h\nabla\nu,\ve{v}_h)\right| \leq  \frac{1}{6}\nu_{\mathrm{min}}\|\nabla\ve{v}_h\|^2 + \frac{3\|\nabla\nu\|_{\infty}^2}{2\nu_{\mathrm{min}}}\|\ve{v}_h\|^2\, ,\\
   &\left|\delta(\nabla q_h,2\nabla^{\top}\ve{v}_h\nabla\nu)\right| \leq   \frac{1}{4}\delta\|\nabla q_h\|^2 + 4\delta\|\nabla\nu\|_{\infty}^2\|\nabla\ve{v}_h\|^2 \, , \\
   &|\delta(\nabla q_h\times\ve{n}, \nu\nabla\times\ve{v}_h)_{\Gamma}| 
\leq \frac{1}{4}\delta\|\nabla q_h \|^2 + 4\delta C\frac{\nu_{\mathrm{max}}^2} {h^2}\|\nabla\ve{v}_h\|^2 \, .
\end{align*}
Hence, we get
\begin{align*}
&A_{\mathrm{GL}}((\ve{v}_h,q_h),(\ve{v}_h,q_h)) \geq  \frac{1}{2}\vvvert{(\ve{v}_h,q_h)}^2 \\ &+ \bigg(1-\frac{3\|\nabla\nu\|^2_{\infty}}{\sigma\nu_{\mathrm{min}}}\bigg)\frac{\sigma}{2}\|\ve{v}_h\|^2 +  \bigg[1- \frac{12\delta}{h^2\nu_{\mathrm{min}}}(h^2\|\nabla\nu\|^2_{\infty} + C\nu_{\mathrm{max}}^2)\bigg]\frac{\nu_{\mathrm{min}}}{3}\|\nabla\ve{v}_h\|^2\, .
\end{align*}
Provided that the reaction coefficient $\sigma$ is large enough \eqref{restriction}, we guarantee coercivity by taking the stabilisation parameter as in \eqref{deltaGL}.}
\begin{rmk}
The restriction \eqref{restriction} is almost identical to the one obtained by \citet{Anaya2021} to prove the coercivity of their vorticity-based formulation. The condition becomes trivial (that is, $\sigma\geq 0$) if $\nu$ is constant. 
\end{rmk}

\section{Convergence}
Due to the consistency of the BVS formulations, Galerkin orthogonality follows easily.
\begin{lemma}[Galerkin orthogonality]
Let $(\ve{u},p)$ be the solution of \eqref{momentum}--\eqref{DirichletBC0}, and\\ $(\ve{u}_h,p_h)\in X_h\times Y_h$ the solution of \eqref{BVSGL}. If $(\ve{u},p)\in H^2(\Omega)\times H^1(\Omega)$, then
\begin{align*}
    A_{\mathrm{GL}}((\ve{u}-\ve{u}_h,p-p_h),(\ve{v}_h,q_h)) = 0 \quad \text{for all} \quad (\ve{v}_h,q_h)\in X_h\times Y_h \, .
\end{align*}
\end{lemma}
\proof{We have
\begin{align*}
    A_{\mathrm{GL}}((\ve{u}-\ve{u}_h,p-p_h),(\ve{v}_h,q_h)) 
    &= A_{\mathrm{GL}}((\ve{u},p),(\ve{v}_h,q_h)) - A_{\mathrm{GL}}((\ve{u}_h,p_h),(\ve{v}_h,q_h))\\
    &= A_{\mathrm{GL}}((\ve{u},p),(\ve{v}_h,q_h)) - (\ve{f},\delta\nabla q_h + \ve{v}_h)\, ,
\end{align*}
where
\begin{align*}
    A_{\mathrm{GL}}((\ve{u},p),(\ve{v}_h,q_h)) 
    &= \sigma(\ve{u},\ve{v}_h)+(\nu\nabla\ve{u},\nabla\ve{v}_h) - (\nabla^{\top}\ve{u}\nabla\nu,\ve{v}_h)  -(p,\nabla\cdot\ve{v}_h) \\
    &+ \delta(\nabla q_h,\nabla p  - 2\nabla^{\top}\ve{u}_h\nabla\nu) + \delta(\nabla q_h\times\ve{n},\nu\nabla\times\ve{u})_{\Gamma}\, .
\end{align*}
Using that $\ve{v}_h\in [H^1_0(\Omega)]^d$, $q_h\in H^1(\Omega)$, $\ve{u}\in [H^1_0(\Omega)\cap H^2(\Omega)]^d$ and $\nabla\cdot\ve{u}=0$, all the integration steps used previously to derive the weak form are valid. Thus, integration by parts yields
\begin{align*}
 \sigma(\ve{u},\ve{v}_h)+(\nu\nabla\ve{u},\nabla\ve{v}_h) - (\nabla^{\top}\ve{u}\nabla\nu,\ve{v}_h)  -(p,\nabla\cdot\ve{v}_h) &= (\sigma\ve{u}-\nabla\cdot(2\nu\nabla^{\mathrm{s}}\ve{u})+\nabla p,\ve{v}_h) \\
 &= (\ve{f},\ve{v}_h)
\end{align*}
and
\begin{align*}
   \delta(\nabla q_h,\nabla p  - 2\nabla^{\top}\ve{u}_h\nabla\nu) + \delta(\nabla q_h\times\ve{n},\nu\nabla\times\ve{u})_{\Gamma} &= \delta(\nabla q_h,\sigma\ve{u}-\nabla\cdot(2\nu\nabla^{\mathrm{s}}\ve{u})+\nabla p)\\
   &=\delta(\nabla q_h,\ve{f})\, ,
\end{align*}
so that 
\begin{align*}
    A_{\mathrm{GL}}((\ve{u}-\ve{u}_h,p-p_h),(\ve{v}_h,q_h)) 
    = (\ve{f},\delta\nabla q_h + \ve{v}_h)\, ,
\end{align*}
completing the proof. The same result holds for the stress-divergence variant $A_{\mathrm{SD}}$.}

With both the Galerkin orthogonality and the stability of our formulations guaranteed, error estimates can be derived using standard techniques. Optimal convergence in the mesh-dependent norm \eqref{norm} can be shown almost exactly as done, for instance, by \citet{John2016} (chapter 4.5.1). Pressure convergence in $H^1(\Omega)$ and $L^2(\Omega)$ can be shown using the techniques by \citet{Araya2023}. 


\section{Numerical experiments}

To assess the BVS method, we conduct simple numerical experiments comparing it to PSPG. We test out two source-free ($\ve{f}=\ve{0}$) analytical solutions of the governing equations in 2D ($d=2$), considering the classical (SD) formulation of the momentum equation (for experiments using the GL version, refer to some of our recent works \cite{Schussnig2021,Barrenechea2024}). The matrix assembly and the  computation of the solutions is done by means of the in-house library MAD (\cite{galarceThesis}, chapter 5), which is based upon the linear algebra library Petsc \cite{petsc}. The  MUltifrontal Massively Parallel sparse direct Solver (MUMPS) is used.

\begin{figure}[!htbp]
    \centering
    \includegraphics[width=0.6\linewidth]{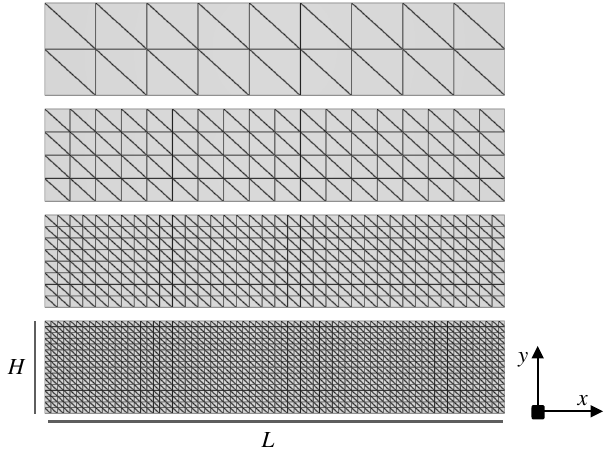}
    \caption{4 coarser structured mesh refinements used for the computational domain $\Omega = (0,5) \times (0,1)$.}
    \label{fig:meshes}
\end{figure}

The two following numerical experiments are computed in a rectangular channel $\Omega = (0,L)\times(0,H)$, $H=1$ and $L=5$, discretised using structured triangular meshes, as depicted in Figure \ref{fig:meshes}. The meshes are refined uniformly to obtain 10 successive configurations. The finest mesh leads to a problem with nearly 3.5 million degrees of freedom.

\section{A reaction-free experiment}
The first numerical test serves to assess accuracy in the particular scenario of a null reaction term. This test highlights the BVS's feature of correcting the spurious boundary condition induced by classical PSPG. For $\sigma=0$ and $\nu(x,y) = a y + b$, we use the following analytic solution:
\begin{equation}
\left\lbrace
\begin{aligned}
\ve{u} &= 
\begin{pmatrix}
    \frac{\kappa}{a} \left[ H - y + \frac{b}{a} \text{ln} \left( \frac{y + b/a}{H + b/a}\right) \right] \\ 0
\end{pmatrix}\, ,\\
p &= \kappa \left( \frac{L}{2} - x \right),
\end{aligned}
\right.
\end{equation}
where $a$, $b$ and $\kappa$ are arbitrary positive parameters, and $H$ and $L$ stand for the height and length of the rectangular domain. For this test, we set $a=b=1$, and $\kappa=0.4$.

\begin{figure}
    \centering
    \subfigure[BVS]{
    \includegraphics[width=0.45\linewidth]{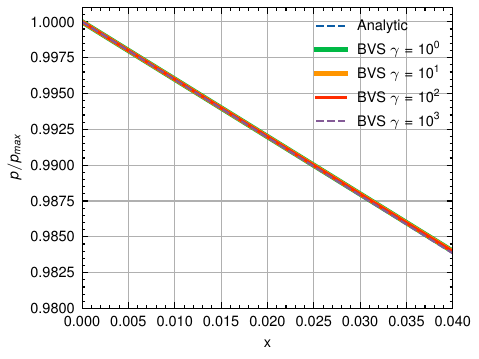}}
    \subfigure[PSPG]{
    \includegraphics[width=0.45\linewidth]{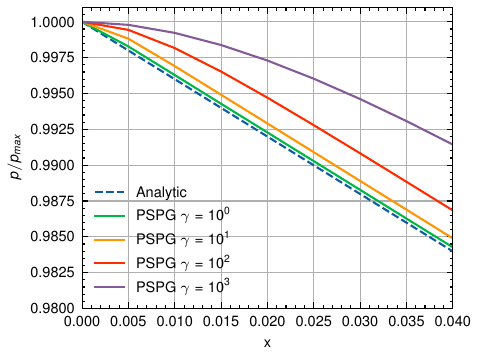}}
    \caption{Centreline pressure for the channel flow with BVS and PSPG. The pressure boundary layer, a numerical artefact, can clearly be seen when using PSPG. The BVS results, on the other hand, show that a careful treatment of the viscous term in the stabilisation leads to a consistent approach, even for over-regularised setups.}
    \label{fig:zoom}
\end{figure}

For this type of problem, the vanishing of the Laplacian in the PSPG stabilisation term induces a spurious boundary condition $\partial_{n} p \approx 0$ on $\Gamma$, as shown by \citet{Schussnig2021} (see their section 6.1). Indeed, the numerical results in Figure \ref{fig:zoom} clearly depict this phenomenon, where we exaggerate its effect by varying the stabilisation term \eqref{BVS} as
\begin{equation}
\delta = \gamma \frac{\nu_{\text{min}}h^2/12 }{h^2 \lVert \nabla{\nu} \rVert_{\infty}^2 + \nu_{\text{max}}^2 }\, ,
\label{eq:stab_gamma}
\end{equation}
with $\gamma > 0$ adjustable. We observe that, even for large $\gamma$, BVS remains consistent and no spurious pressure boundary condition is observed. For PSPG, on the other hand, increasing $\gamma$ clearly expands the numerical boundary layer.

\section{Generalised Stokes experiment}
\begin{figure}
    \centering
    \includegraphics[width=0.7\linewidth]{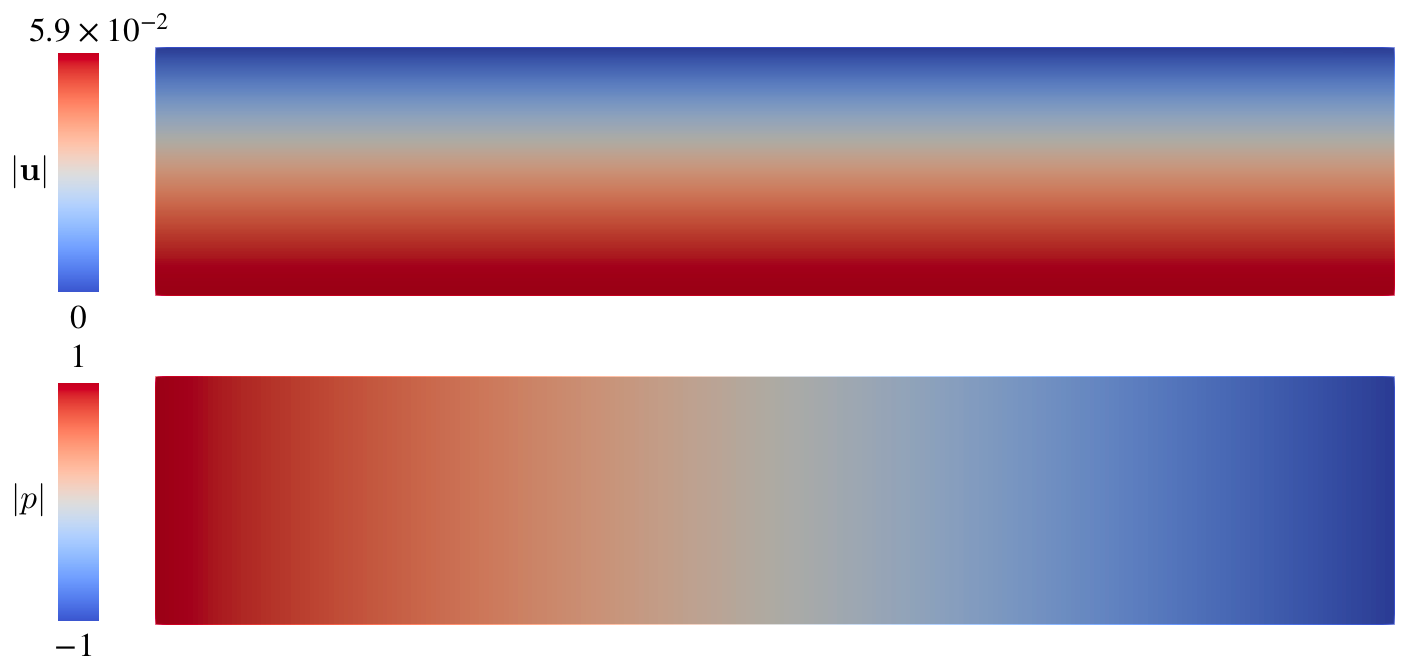}
    \caption{Numerical solution using BVS for the generalised Stokes experiment}
    \label{fig:sol_exp2}
\end{figure}

\begin{figure}[!htbp]
    \centering
    \includegraphics[width=\linewidth]{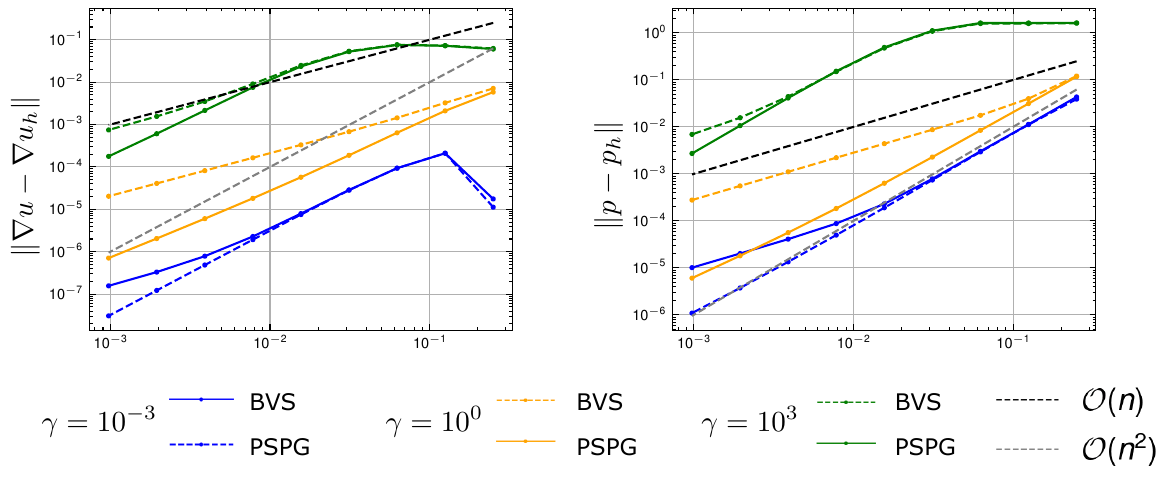}
    \caption{Solution convergence for $\ve{u}$ and $p$ in generalized Stokes experiment and its sensitivity with respect to the parameter $\gamma$ in \eqref{eq:stab_gamma}.}
    \label{fig:convergence_linear}
\end{figure}

We now consider another channel flow solution, derived for $\sigma=2$ and $\nu(x,y)=(y+b)^2$:
\begin{equation}
\left\lbrace
\begin{aligned}
\ve{u}  &= 
    \begin{pmatrix}
       \frac{\kappa}{2} \left[ 1 - \frac{2(y+b)}{\beta b^3} - \frac{1}{\beta (y+b)^2} \right] \\ 0
    \end{pmatrix}\\
    p &= \kappa \left( \frac{L}{2} - x \right)
\end{aligned}
\right. ,
\end{equation}
where
\begin{align*}
    \beta = \frac{2(H+b)}{b^3} + \frac{1}{(H+b)^2}.
\end{align*}
The parameters are set so as to keep the values of $p$ within $[-1,1]$: $\kappa = 0.4$, $b = 1$, again with $L = 5$ and $H = 1$. This leads to the solutions depicted in Figure \ref{fig:sol_exp2}. We assess the finite element convergence numerically in Figure \ref{fig:convergence_linear}, using different stabilisation values \eqref{eq:stab_gamma}. As expected from classical theory, both methods converge at least linearly. For PSPG, however, some pressure super-convergence is observed, which is also a well-known phenomenon for very regular solutions and meshes \cite{Eichel2011}. Yet, this is highly problem-dependent, and in other cases BVS may yield better convergence, see e.g.~some results we previously reported \cite{Schussnig2021}. In any case, such super-convergence cannot be expected for more realistic applications with less smooth solutions, where the initial super-convergence usually breaks down after some refinement \cite{Pacheco2022CAMWA}. 

\begin{figure}
    \centering
    \subfigure[BVS]{
    \includegraphics[width=0.45\linewidth]{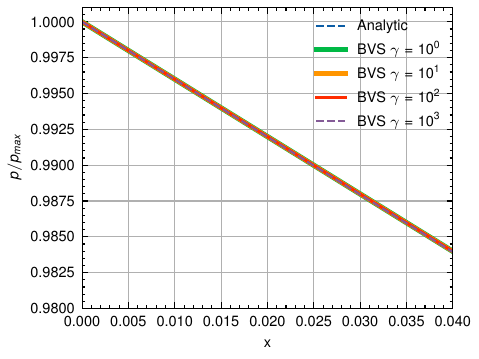}}
    \subfigure[PSPG]{
    \includegraphics[width=0.45\linewidth]{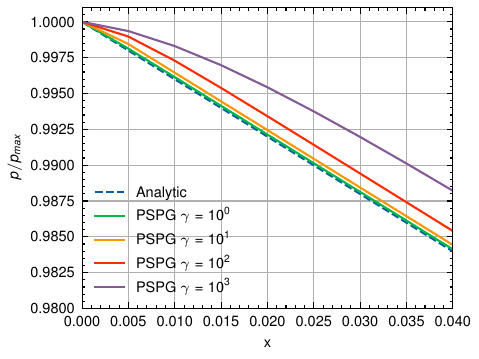}}
    \caption{Spurious pressure boundary layers induced by PSGP and corrected by BVS.}
    \label{fig:bl_ex2}
\end{figure}

Finally, we depict again the boundary layer phenomenon for PSPG in Figure \ref{fig:bl_ex2} and verify (as with Example 1) how BVS is able to correct this issue thanks to its simple boundary integral.

\newpage
\section{Concluding remarks}
In this article, we have analysed the well-posedness of a pressure stabilisation method that retains full consistency even for the (continuous) lowest-order pressure-velocity finite element pairs. Through the analysis, we have derived an expression for the stabilisation parameter, which can be seen as a variable-viscosity generalisation of the well-known PSPG parameter. Our numerical methods show that, although our formulation -- as expected -- does not improve the convergence order in comparison to PSPG, it completely eliminates the issue of spurious pressure boundary layers. This is particularly relevant in applications where the pressure on the boundary is a key quantity, such as in aero- and hemodynamics.

\section*{Acknowledgements}
 The first author acknowledges the research funding VINCI PUCV DI-Iniciación 039.482/2024. The corresponding author acknowledges funding by the Federal Ministry of Education and Research (BMBF) and the Ministry of Culture and Science of the German State of North Rhine-Westphalia (MKW) under the Excellence Strategy of the Federal Government and the Länder.



\end{document}